\documentstyle[a4]{article}

\begin{document}

\title{On the Cauchy- and periodic boundary value problem for a certain class of derivative nonlinear Schr\"odinger equations}

\author{Axel Gr\"unrock \\ Bergische Universit\"at - Gesamthochschule Wuppertal \\ Gau{\ss}stra{\ss}e 20 \\ D-42097 Wuppertal \\ Germany \\ e-mail Axel.Gruenrock@math.uni-wuppertal.de}

\date{}

\maketitle

\begin{abstract}
The Cauchy- and periodic boundary value problem for the nonlinear Schr\"o- dinger equations in $n$ space dimensions
\[u_t - i\Delta u = (\nabla \overline{u})^{\beta}, \hspace{1cm}|\beta|=m \ge 2, \hspace{1cm}u(0)=u_0 \in H^{s+1}_x\]
is shown to be locally well posed for $s > s_c := \frac{n}{2} - \frac{1}{m-1}$ , $s \ge 0$. In the special case of space dimension $n=1$ a global $L^2$-result is obtained for NLS with the nonlinearity $N(u)= \partial_x (\overline{u} ^2)$. The proof uses the Fourier restriction norm method.
\end{abstract}

\newcommand{\R}{\mbox{${\bf R}$}}
\newcommand{\N}{\mbox{${\bf N}$}}
\newcommand{\C}{\mbox{${\bf C}$}}
\newcommand{\T}{\mbox{${\bf T}$}}
\newcommand{\Z}{\mbox{${\bf Z}$}}
\newcommand{\r}{\mbox{${\longrightarrow}$}}
\newcommand{\imp}{\mbox{${\Rightarrow}$}}
\newcommand{\mip}{\mbox{${\Leftarrow}$}}
\newcommand{\iso}{\mbox{${\stackrel{\sim}{\longrightarrow}}$}}
\newcommand{\F}{\mbox{${\cal F}$}}
\newcommand{\Fx}{\mbox{${\cal F}_x$}}
\newcommand{\Ft}{\mbox{${\cal F}_t$}}
\newcommand{\Sn}{\mbox{${\cal S} (\R^{n+1})$}}
\newcommand{\So}{\mbox{${\cal S} (\R)$}}
\newcommand{\Y}{\mbox{${\cal Y} (\R \times \T^n)$}}
\newcommand{\U}{\mbox{$ U_{\phi}$}}
\newcommand{\x}{\mbox{$ X^+_{s,b} $}}
\newcommand{\y}{\mbox{$ Y_{s}$}}
\newcommand{\xm}{\mbox{$ X^-_{s,b} $}}
\newcommand{\ym}{\mbox{$ Y^-_{s}$}}
\newcommand{\xpm}{\mbox{$ X^{\pm}_{s,b} $}}
\newcommand{\ypm}{\mbox{$ Y^{\pm}_{s}$}}
\newcommand{\yy}[1]{\mbox{$ Y_{#1} $}}
\newcommand{\X}[1]{\mbox{$ X_{s,#1}$}}
\newcommand{\Xpm}[1]{\mbox{$ {X^{\pm}}_{s,#1}$}}
\newcommand{\XX}[2]{\mbox{$ X^+_{#1,#2} $}}
\newcommand{\XXm}[2]{\mbox{$ X^-_{#1,#2} $}}
\newcommand{\XXpm}[2]{\mbox{$ X^{\pm}_{#1,#2} $}}
\newcommand{\XXX}[3]{\mbox{$ X_{#1,#2}(#3) $}}
\newcommand{\nfx}{\mbox{$ \| f \| _{ X_{s,b}} $}}
\newcommand{\qfx}{\mbox{$ {\| f \|}^2_{ X_{s,b}} $}}
\newcommand{\nf}[1]{\mbox{$ \| f \| _{ #1} $}}
\newcommand{\qf}[1]{\mbox{$ {\| f \|}^2_{ #1} $}}
\newcommand{\nx}[1]{\mbox{$ \| #1 \| _{ X_{s,b}} $}}
\newcommand{\qx}[1]{\mbox{$ {\| #1 \|}^2_{X_{s,b} } $}}
\newcommand{\n}[2]{\mbox{$ \| #1 \| _{ #2} $}}
\newcommand{\nk}[3]{\mbox{$ \| #1 \| _{ #2} ^{#3}$}}
\newcommand{\q}[2]{\mbox{$ {\| #1 \|}^2_{#2} $}}
\newcommand{\conv}[2]{\mbox{${\mbox{\Huge{$\ast$}}}_{_{\!\!\!\!\!\!\!\!\!{#1}}}^{^{\!\!\!\!\!\!{#2}}}$}}

\pagestyle{plain}
\rule{\textwidth}{0.5pt}

\newtheorem{lemma}{Lemma}[section]
\newtheorem{bem}{Remark}[section]
\newtheorem{bsp}{Example}[section]
\newtheorem{definition}{Definition}[section]
\newtheorem{kor}{Corollary}[section]
\newtheorem{satz}{Theorem}[section]
\newtheorem{prop}{Proposition}[section]

\section{Introduction and main results}

In this paper we prove local (in time) wellposedness of the initial value and periodic boundary value problem for the following class of derivative nonlinear Schr\"odinger equations
\[u_t - i\Delta u = (\nabla \overline{u})^{\beta}, \hspace{2cm}u(0)=u_0 \in H^{s+1}_x .\]
Here the initial value $u_0$ belongs to the Sobolev space $H^{s+1}_x =H^{s+1}_x (\R ^n)$ or $H^{s+1}_x =H^{s+1}_x (\T ^n)$, $\beta \in \N_0^n$ is a multiindex of length $|\beta|=m \geq 2$ and we can admit all values of $s$ satisfying
\[s > s_c := \frac{n}{2} - \frac{1}{m-1}, \hspace{3cm} s \ge 0 .\]
The same arguments give local wellposedness for the problem
\[u_t - i\Delta u = \partial_j (\overline{u}^{m}), \hspace{2cm}u(0)=u_0 \in H^{s}_x\]
with the same restrictions on $s$ as above. In the special case of a quadratic nonlinearity in one space dimension (i. e. $m=2,\,\,n=1$) we can reach the value $s=0$. Employing the conservation of $\n{u(t)}{L^2_x}$ in this case, we obtain global wellposedness for
\[u_t - i \partial_x^2 u = \partial_x (\overline{u}^{2}), \hspace{2cm}u(0)=u_0 \in L^{2}_x\]
respectively for
\[u_t - i  \partial_x^2 u= (\partial_x \overline{u})^{2}, \hspace{2cm}u(0)=u_0 \in \dot{H}^{1}_x .\]
To prove our results, we use the Fourier restriction norm method as it was introduced in \cite{B93} and further developped in \cite{KPV96} and \cite{GTV97}. In particular, we will use the function spaces $\xpm = \exp{(\pm it\Delta)} H^b_t(H^s_x)$ equipped with the norms
\begin{eqnarray*}
\n{f}{\xpm}= \n{\exp{(\mp it\Delta)}f}{H^b_t(H^s_x)}=\n{<\xi>^s <\tau \pm |\xi|^2>^b \F f}{L^2_{\xi \tau}} \\
= (\int \mu (d \xi) d \tau <\xi>^{2s}<\tau \pm |\xi|^2>^{2b}|\F f (\xi,\tau)|^2)^{\frac{1}{2}} \hspace{1,5cm}
\end{eqnarray*}
as well as the auxiliary norm
\begin{eqnarray*}
\n{f}{\y}=\n{<\xi>^s <\tau + |\xi|^2>^{-1} \F f}{L^2_{\xi }(L^1_{\tau})} \hspace{1cm} \\
= (\int \mu (d \xi)<\xi>^{2s} (\int d \tau <\tau + |\xi|^2>^{-1}|\F f (\xi,\tau)|)^2)^{\frac{1}{2}}
\end{eqnarray*}
introduced in \cite{GTV97}. Here $\F$ denotes the Fourier transform in space and time, $\mu$ is the Lebesgue measure on $\R^n$ in the nonperiodic respectively the counting measure on $\Z^n$ in the periodic case, and we use the notation $<x>=(1+|x|^2)^{\frac{1}{2}}$. Observe that $\n{\overline{f}}{\x}=\n{f}{\xm}$. To give a precise formulation of our results, we also need the restriction norm spaces $\xpm (I) = \exp{(\pm it\Delta)} H^b_t(I, H^s_x)$ with norms
\[\n{f}{\xpm (I)}= \inf \{\n{\tilde{f}}{\xpm} : \tilde{f} \in \xpm \,\,\, ,\,\,\,\tilde{f}|_I=f\} .\]
Now our results read as follows:
\begin{satz}Let $n=1$ and $s\ge 0$. Then there exists a $T=T(\n{u_0}{H^s_x})>0$, so that there is a unique solution $u \in \x ([-T,T])$ of the initial value (periodic boundary value) problem
\[u_t - i \partial_x^2 u = \partial_x (\overline{u}^{2}) , \hspace{2cm}u(0)=u_0 \in H^{s}_x .\]
This solution satisfies $u \in C^0([-T,T],H^s_x)$ and for any $T' < T$ the mapping \\ $f: H^s_x \r \x ([-T',T'])\,,\,u_0 \mapsto u$ (Data upon solution) is locally Lipschitz continuous. For $s=0$, by the conservation of $\n{u(t)}{L^2_x}$, this solution extends globally, in this case it is $u \in C^0(\R,L^2_x)$.
\end{satz}
Using the isomorphism $\partial_x: \dot{H}^{1}_x \r L^2_x$ we obtain the following 
\begin{kor} The Cauchy- and the periodic boundary value problem
\[u_t - i \partial_x^2 u = (\partial_x \overline{u})^{2} , \hspace{2cm}u(0)=u_0 \in \dot{H}^{1}_x\]
is globally well posed.
\end{kor}
Our next result is much more general:

\begin{satz}Let $m,n \in \N$, $m\ge 2$ and $m+n \ge 4$. Then for $s>s_c$ there exists a $T=T(\n{u_0}{H^s_x})>0$ and a unique solution $u \in \x ([-T,T])$ of the initial value (periodic boundary value) problem
\[u_t - i\Delta u = \partial_j (\overline{u}^{m}) , \hspace{2cm}u(0)=u_0 \in H^{s}_x .\]
This solution is persistent and for any $T' < T$ the mapping Data upon solution from $H^s_x$ to $\x ([-T',T'])$ is locally Lipschitz continuous. \\ For any $\beta \in \N_0^n$ with $|\beta|= m$ and under the same assumptions on $m,n,s$   the Cauchy problem and the periodic boundary value problem
\[u_t - i\Delta u = (\nabla \overline{u})^{\beta} , \hspace{2cm}u(0)=u_0 \in H^{s+1}_x\]
is locally well posed in the same sense.
\end{satz}

$Remarks:$ 1. The special case in Theorem 1.2, where $n=1$, $m=3$ and $s>0$, has already been proved for the nonperiodic case by H. Takaoka, see Thm. 1.2 in \cite{T}.

2. A standard scaling argument suggests, that our result is optimal as long as we are not dealing with the critical case $s=s_c$. In fact, if $u$ is a solution of the first problem in Theorem 1.2 with initial value $u_0 \in H^s_x(\R^n)$, then so is $u_{\lambda}$, defined by $u_{\lambda}(x,t)=\lambda^{\frac{1}{m-1}}u(\lambda x,\lambda ^2 t)$, with initial value  $u_{\lambda}^0(x)= u_0(\lambda x)$, and $\n{u_{\lambda}^0}{\dot{H}^{s_c}_x}(\R^n)$ is independent of $\lambda$.

\vspace{0,5cm}

As mentioned above, we use the Fourier restriction norm method. We shall assume this method to be known in order to concentrate on the derivation of the crucial nonlinear - in our case in fact multilinear - estimates, see Theorem 3.1 (3.2) below, corresponding to Theorem 1.1 (1.2) of this introduction. Our proofs rely heavily on the following interpolation property of the $\xpm$-spaces: We have
\[(\XXpm{s_0}{b_0},\XXpm{s_1}{b_1})_{[\theta]} = \xpm \,\,\,,\]
whenever for $\theta \in [0,1]$ it holds that $s=(1-\theta )s_0 + \theta s_1$, $b=(1-\theta )b_0 + \theta b_1$. Here $[\theta]$ denotes the complex interpolation method. Moreover we will make extensive use of the fact, that with respect to the inner product on $L^2_{xt}$ the dual space of $\xpm$ is given by $\XXpm{-s}{-b}$. Another useful tool to derive the nonlinear estimates is the following

\begin{lemma}\label{l11}
For $1 \leq i \leq m$ let $u_i \in \XXpm{s}{b_i}$ and 
\[f_i(\xi, \tau) := <\tau \pm |\xi|^2>^{b_i}<\xi>^{s}\F u_i(\xi, \tau).\] Then with
$ d\nu := \mu (d\xi _1 .. d\xi _{m-1}) d \tau _1 .. d \tau _{m-1}$ and $ \xi = \sum_{i=1}^m \xi _i, \,\,\tau = \sum_{i=1}^m \tau _i$ the following identities are valid:
\[\F D^{\beta}(\prod_{i=1}^m D^{\beta _i}u_i)(\xi,\tau)=c \xi ^{\beta} \int d \nu \prod_{i=1}^m \xi_i^{\beta_i}<\tau_i \pm |\xi_i|^2>^{-b_i}<\xi_i>^{-s}f_i(\xi_i,\tau_i)\]
as well as
\begin{itemize}
\item[a)] $\n{D^{\beta}(\prod_{i=1}^m D^{\beta _i}u_i)}{\Xpm{b'}}=$ 
\end{itemize}
$c\n{<\tau \pm |\xi|^2>^{b'}<\xi>^s\xi ^{\beta} \int d \nu \prod_{i=1}^m \xi_i^{\beta_i}<\tau_i \pm |\xi_i|^2>^{-b_i}<\xi_i>^{-s}f_i(\xi_i,\tau_i) }{L^2_{\xi, \tau}}$ and
\begin{itemize}
\item[b)] $\n{D^{\beta}(\prod_{i=1}^m D^{\beta _i}u_i)}{\y}=$ 
\end{itemize}
 $c\n{<\tau + |\xi|^2>^{-1}<\xi>^s\xi ^{\beta} \int d \nu \prod_{i=1}^m \xi_i^{\beta_i}<\tau_i \pm |\xi_i|^2>^{-b_i}<\xi_i>^{-s}f_i(\xi_i,\tau_i) }{L^2_{\xi}(L^1_{ \tau})}$.
\end{lemma}

Proof: For the convolution of $m$ functions $g_i,\,\,1 \leq i \leq m$, we have with \\ $x = \sum_{i=1}^m x_i$
\[\conv{i=1}{m} g_i (x) = \int \mu(d x_1 .. dx_{m-1}) \prod_{i=1}^m g_i (x_i)\]
Hence by the properties of the Fourier transform the following holds true with \\ $\xi = \sum_{i=1}^m \xi _i, \,\,\tau = \sum_{i=1}^m \tau _i$:
\begin{eqnarray*}
& &\F D^{\beta}(\prod_{i=1}^m D^{\beta _i}u_i)(\xi,\tau)\\
&=& c \xi ^{\beta} (\conv{i=1}{m} \xi ^{\beta _i} \F u_i)(\xi,\tau) \\
&=& c \xi ^{\beta} (\conv{i=1}{m} \xi ^{\beta _i}<\tau \pm |\xi|^2>^{-b_i}<\xi>^{-s}f_i)(\xi,\tau) \\
&=& c \xi ^{\beta} \int d \nu \prod_{i=1}^m \xi_i^{\beta_i}<\tau_i \pm |\xi_i|^2>^{-b_i}<\xi_i>^{-s}f_i(\xi_i,\tau_i) .
\end{eqnarray*}
From this we obtain a) because of 
\[\n{D^{\beta}(\prod_{i=1}^m D^{\beta _i}u_i)}{\XXpm{s}{b'}}= \n{<\tau \pm |\xi|^2>^{b'}<\xi>^s \F D^{\beta}(\prod_{i=1}^m D^{\beta _i}u_i)}{L^2_{\xi, \tau}}\]
and b) because of
\[\n{D^{\beta}(\prod_{i=1}^m D^{\beta _i}u_i)}{\y}= \n{<\tau + |\xi|^2>^{-1}<\xi>^s \F D^{\beta}(\prod_{i=1}^m D^{\beta _i}u_i)}{L^2_{\xi}(L^1_{\tau})}.\]
$\hfill \Box$

$Remark:$ The previous Lemma has some simple but important consequences: First of all it shows, that the estimate
\[\n{D^{\beta}(\prod_{i=1}^m D^{\beta _i}u_i)}{\XXpm{s}{b'}} \leq c \prod_{i=1}^m \n{u_i}{\XXpm{s}{b_i}}\]
holds true, iff
\begin{eqnarray*}
\n{<\tau \pm |\xi|^2>^{b'}<\xi>^s \xi^{\beta} \int d \nu \prod_{i=1}^m \xi_i^{\beta_i}<\tau_i \pm |\xi_i|^2)>^{-b_i}<\xi_i>^{-s}f_i(\xi_i,\tau_i)}{L^2_{\xi, \tau}} 
 \nonumber \\
\leq c \prod_{i=1}^m \n{f_i}{L^2_{\xi, \tau}} .\hspace{5cm}
\end{eqnarray*}
In order to prove the latter one may assume without loss of generality, that \\ $\xi^{\beta}\prod_{i=1}^m \xi_i^{\beta_i}f_i(\xi_i,\tau_i) \geq 0$. (This will be used throughout section 3.) Because of 
\[<\xi>=<\sum_{i=1}^m \xi _i> \leq \sum_{i=1}^m <\xi _i>\leq \prod_{i=1}^m <\xi _i>\]
it follows, that if the above estimates hold for some $s \in \R$, then they are valid for every $s' \geq s$. Corresponding statements are true with $\XXpm{s}{b'}$ replaced by $\y$.

\vspace{0,5cm}

In order to extract a positive power of $T$ (the lifespan of the local solutions) from the nonlinear estimates the following Lemma will be useful:
\begin{lemma}\label{l12} Let $\psi \in C^{\infty}_0 (\R)$ be a smooth characteristic function of $[-1,1]$ and $\psi_{\delta}(t):= \psi(\frac{t}{\delta})$ for $0< \delta \leq 1$. Then, if $1/2 >b>b' \geq 0$ or $0 \geq b>b'> -1/2$, the following estimate holds:
\[\n{\psi_{\delta}f}{\XXpm{s}{b'}} \leq c \delta^{b-b'}\n{f}{\xpm}\,\, .
\]
\end{lemma}

Proof: First assume $1/2 >b>b' \geq 0$. Then for $g \in H^b_t$ we use the generalized Leibniz rule and Sobolev's embedding theorem with $b-b' = 1/p$ and $b= 1/q$ as well as $1/p +1/p' = 1/q +1/q' = 1/2$ to obtain:
\begin{eqnarray*}
\n{\psi_{\delta}g}{H^{b'}_t} &\leq &c (\n{\psi_{\delta}}{L^p_t} \n{g}{H^{b',p'}_t}+\n{\psi_{\delta}}{H^{b',q}_t}\n{g}{L^{q'}_t}) \\
&\leq& c(\n{\psi_{\delta}}{L^p_t}+\n{\psi_{\delta}}{H^{b',q}_t})\n{g}{H^b_t}\\
&\leq& c\delta^{b-b'}\n{g}{H^b_t}\,\,.
\end{eqnarray*}
From this for $f\in \xpm$ it follows with $U(t)= \exp{(it\Delta)}$:
\begin{eqnarray*}
\n{\psi_{\delta}f}{\xpm} &= &\n{U (\mp .)\psi_{\delta}f}{H^b_t(\R,H^s_x)}\\
&= &\n{\psi_{\delta} U (\mp .)f}{H^b_t(\R,H^s_x)}\\
&\leq &c\delta^{b-b'}\n{U (\mp .)f}{H^b_t(\R,H^s_x)}=c \delta^{b-b'}\nf{\xpm}\,\,.
\end{eqnarray*}
Exchanging $b$ and $b'$ the case $0 \geq b>b'> -1/2$ follows from this by duality.
$\hfill \Box$

\vspace{0,5cm}

An important ingredient of our proofs are the Strichartz type estimates for the homogeneous linear Schr\"odinger equation. Here we will make use of some of the results derived in \cite{B93}, in a slightly modified version, which we prove in section 2. This will enable us to give proofs of the nonlinear estimates that cover both, the periodic and the nonperiodic case. To obtain the nonlinear estimates necessary for Thm. 1.2 (see Thm. 3.2 below) we follow the ideas of section 5 in \cite{B93} - essentially we present a simplified version of the proof given there. To do so, we took some instructive hints from \cite{G96}, section 5. In particular, we do use Hilbert space norms instead of Besov-type norms as in \cite{B93}. Perhaps it is worthwile to mention, that for the nonperiodic case there is a much easier proof of Theorem 3.2 below, using the full strength of the Strichartz estimates in this case.

\vspace{0,5cm}

{\bf{Acknowledgement:}} I want to thank Professor Hartmut Pecher for numerous helpful conversations.

\section{Some Strichartz type estimates for the Schr\"odinger equation in the periodic case}

In this section we prove some of the Strichartz type estimates for the Schr\"odinger equation in the periodic case, which were discoverd by Bourgain in \cite{B93}. All the following estimates are essentially contained in \cite{B93}. Since we want to use them in the form of an embedding of the type $L^p_t(L^q_x) \subset \x$, where we have spaces of functions being periodic in the space- but not in the time-variable, we shall give new proofs for these estimates, using some of the arguments from \cite{B93} as well as from \cite{KPV96} and \cite{GTV97}.

\begin{lemma}\label{l21} (cf. \cite{B93}, Prop. 2.6) Let $n=1$. Then for any $b> \frac{3}{8}$ and for any $b' <  -\frac{3}{8}$ the following estimates hold:
\begin{itemize}
\item[i)]$\nf{L^4_{xt}} \leq c \nf{\XX{0}{b}}$
\item[ii)]$\nf{\XX{0}{b'}} \leq c \nf{L^{\frac{4}{3}}_{xt}}$
\end{itemize}
\end{lemma}

Proof (cf. \cite{KPV96}, Lemma 5.3): Clearly, ii) follows from i) by duality. To see i), we shall show first, that
\[ \sup_{(\xi, \tau) \in \Z \times \R} S(\xi, \tau) < \infty\]
for
\begin{eqnarray*}
S(\xi, \tau) & = & \sum_{\xi_1 \in \Z}<\tau + \xi_1^2 + (\xi -\xi_1)^2 >^{1 - 4b} \\
& \leq & c \sum_{\xi_1 \in \Z}<4 \tau + (2 \xi_1)^2 + (2 (\xi -\xi_1))^2 >^{1 - 4b}.
\end{eqnarray*}
With $k= 2 \xi_1 - \xi \in \Z$ we have
\[k+\xi=2 \xi_1 , \,\,k-\xi=2 (\xi_1 -\xi) \,\,\mbox{and}\,\,(2 \xi_1)^2 + (2 (\xi -\xi_1))^2 = 2(\xi ^2 + k^2),\]
hence 
\begin{eqnarray*}
S(\xi, \tau) & \leq & c \sum_{k \in \Z}<4 \tau + 2 \xi^2 + 2 k^2 >^{1 - 4b} \\
& \leq & c \sum_{k \in \Z}<  k^2 - |2 \tau +  \xi^2|>^{1 - 4b} \\
& \leq & c \sum_{k \in \Z}<  (k - x_0)(k + x_0)>^{1 - 4b},
\end{eqnarray*}
where $x_0^2 = |2 \tau +  \xi^2|$. Now there are at most four numbers $k \in \Z$ with $|k - x_0|<1$ or $|k + x_0|<1$. For all the others we have
\[<k - x_0><k + x_0> \leq c <(k - x_0)(k + x_0)> .\]
Cauchy-Schwarz' inequality gives
\begin{eqnarray*}
S(\xi, \tau) & \leq & c + c \sum_{k \in \Z}(<k - x_0><k + x_0>)^{1 - 4b} \\
& \leq & c + c (\sum_{k \in \Z}<k - x_0>^{2(1 - 4b)})^{\frac{1}{2}}(\sum_{k \in \Z}<k + x_0>^{2(1 - 4b)})^{\frac{1}{2}} \leq c \,\,\,,
\end{eqnarray*}
provided $2(1 - 4b) < -1$, that is $b>\frac{3}{8}$. Without loss of generality we can assume $b \in (\frac{3}{8},\frac{1}{2})$. Then with Lemma 4.2 in \cite{GTV97} we obtain from the above, that
\[ \sup_{(\xi, \tau) \in \Z \times \R}\sum_{\xi_1 \in \Z} \int d \tau_1 <\tau_1 + \xi_1^2 >^{-2b}<\tau - \tau_1 + (\xi -\xi_1)^2 >^{-2b} < \infty \,\,\,.\]
Using Cauchy-Schwarz' inequality and Fubini's theorem as in \cite{KPV96} (there: proof of Theorem 2.2) we arrive at
\begin{eqnarray*}
\n{\sum_{\xi_1 \in \Z} \int d \tau_1 <\tau_1 + \xi_1^2>^{-b} f(\xi_1, \tau_1)<\tau - \tau_1 + (\xi - \xi_1)^2>^{-b} g(\xi - \xi_1,\tau - \tau_1)}{L^2_{\xi, \tau}} \\
\leq c \n{f}{L^2_{\xi, \tau}}\n{g}{L^2_{\xi, \tau}} .\hspace{5cm}
\end{eqnarray*}
Now by Lemma \ref{l11} from the introduction it follows that
\[\n{u_1 u_2}{L^2_{xt}} \leq c \n{u_1}{\XX{0}{b}}\n{u_2}{\XX{0}{b}}\,\,.\]
Taking $u_1 =u_2 =u$, we get
\[\n{u}{L^4_{xt}} = \n{u^2}{L^2_{xt}}^{\frac{1}{2}} \leq c\n{u}{\XX{0}{b}}\,\,.\]
$\hfill \Box$

In the sequel we shall make use of the following number theoretic results concerning the number of solutions of certain Diophantine equations:

\begin{prop}\label{p21}
\begin{itemize}
\item[i)] For all $\epsilon > 0$ there exists a constant $c=c(\epsilon)$ with 

$a(r,3):=\#\{(k_1,k_2)\in\Z^2:3k_1^2+k_2^2=r\in \N\}\leq c <r>^{\epsilon}$.
\item[ii)] For all $\epsilon > 0$ there exists a constant $c=c(\epsilon)$ with 

$a(r,1):= \# \{(k_1,k_2) \in \Z^2 :  k_1^2 + k_2^2 = r \in \N\} \leq c <r>^{\epsilon}$.
\item[iii)] Let $n \geq 3$. Then for all $\epsilon > 0$ there exists a constant $c=c(\epsilon)$ with 

$ \# \{k \in \Z^n  : |k|^2 = r \in \N\} \leq c <r>^{\frac{n-2}{2} + \epsilon}$.
\end{itemize}
\end{prop}

Quotation/Proof: i) $a(r,3)$ is calculated explicitly in \cite{P}, Satz 6.2: It is
\[a(r,3)=2 (-1)^r \sum_{d|r} (\frac{d}{3})\,\,.\]
Here $(\frac{d}{p})$ denotes the Legendre-symbol taking values only in $\{0, \pm 1 \}$. Thus $a(r,3)$ can be estimated by the number of divisors of $r$, which is bounded by $c <r>^{\epsilon}$, see \cite{HW}, Satz 315. For ii), see Satz 338 in \cite{HW}. iii) follows from ii) by induction, writing $\{k \in \Z^n  : |k|^2 = r \in \N\} = \bigcup_{k_n^2 \leq r} \{(k',k_n)  : |k'|^2 = r-k_n^2 \}$.

\vspace{0,5cm}

The following Lemma corresponds to Prop. 2.36 in \cite{B93}:

\begin{lemma}\label{l22} Let $n=1$. Then for all $s>0$ and $b>\frac{1}{2}$ there exists a constant $c = c(s,b)$, so that the following estimate holds:
\[\nf{L^6_{xt}} \leq c \nf{\XX{s}{b}}\,\,.\]
\end{lemma}

Proof: As in the proof of the previous lemma, we start by showing that
\[ \sup_{(\xi, \tau) \in \Z \times \R} S(\xi, \tau) < \infty\,\,,\]
where now (with $\xi_3 = \xi -\xi_1 - \xi_2$)
\begin{eqnarray*}
S(\xi, \tau) & = & \sum_{\xi_1 , \xi_2 \in \Z}<\tau + \xi_1^2 + \xi_2^2 + \xi_3^2 >^{-2b} <\xi_1>^{-2s}<\xi_2>^{-2s}<\xi_3>^{-2s}\\
& \leq & \!\!\!c \!\!\!\sum_{\xi_1 , \xi_2 \in \Z}\!\!\!<\!\!9 \tau + (3\xi_1)^2 + (3\xi_2)^2 + (3\xi_3)^2 \!\!>^{-2b} <\!\!(3\xi_1)^2 + (3\xi_2)^2 + (3\xi_3)^2\!\!>^{-s}.
\end{eqnarray*}
Taking $k_1 = 3(\xi_1 + \xi_2) - 2\xi$ and $k_2 = 3(\xi_1 - \xi_2)$ as new indices, we have
\[3\xi_1 = \frac{1}{2}(k_1 + k_2) + \xi, \,\,\,3\xi_2 = \frac{1}{2}(k_1 - k_2) + \xi\,\,\,\mbox{and}\,\,\,3 \xi_3 = \xi -k_1 .\]
From this we get
\[(3\xi_1)^2 + (3\xi_2)^2 + (3\xi_3 )^2 = \frac{1}{2}(3 k_1^2 + k_2^2) + 3\xi^2 .\]
It follows
\begin{eqnarray*}
S(\xi, \tau) & \leq & c \sum_{k_1 , k_2 \in \Z}<9 \tau + 3\xi^2 + \frac{1}{2}(3 k_1^2 + k_2^2)>^{-2b} <\frac{1}{2}(3 k_1^2 + k_2^2)>^{-s} \\
& \leq & c \sum_{r \in \N_0} \sum_{3 k_1^2 + k_2^2=r} <9 \tau + 3\xi^2 + \frac{r}{2}>^{-2b} <\frac{r}{2}>^{-s} \\
& \leq & c \sum_{r \in \N_0} <9 \tau + 3\xi^2 + \frac{r}{2}>^{-2b},
\end{eqnarray*}
where in the last step we have used part i) of the above proposition. Since we have demanded $b>\frac{1}{2}$, the introducing claim follows. Again, we use Lemma 4.2 from \cite{GTV97} to obtain
\[\sup_{(\xi, \tau) \in \Z \times \R} \int d \nu \prod_{i=1}^3 <\tau_i + \xi_i^2>^{-2b}<\xi_i>^{-2s} < \infty\]
with $\int d \nu = \int d \tau _1 d \tau _2 \sum_{\xi_1 , \xi_2 \in \Z}$ and $(\tau, \xi)= \sum_{i=1}^3 (\tau_i, \xi_i)$.
Now Cauchy-Schwarz and Fubini are applied to obtain
\[\n{\int d \nu \prod_{i=1}^3 <\tau_i + \xi_i^2>^{-b}<\xi_i>^{-s}f_i(\xi_i,\tau_i)}{L^2_{\xi, \tau}} 
\leq c \prod_{i=1}^3 \n{f_i}{L^2_{\xi, \tau}}\,\,\,. \]
Lemma \ref{l11} gives
\[\n{\prod_{i=1}^3 u_i}{L^2_{xt}} \leq c \prod_{i=1}^3 \n{u_i}{\x}\,\,\,.\]
Because of $\| u \|^3_{L^6_{xt}} = \| u^3 \|_{L^2_{xt}}$ the proof is complete.
$\hfill \Box$

\begin{kor}\label{k21} Let $n=1$:
\begin{itemize}
\item[a)] For all H\"older- and Sobolevexponents $p,\,\,q,\,\,s$ and $b$ satisfying
\[0 \leq \frac{1}{p} \leq \frac{1}{6},\,\,\,0 < \frac{1}{q} \leq \frac{1}{2} - \frac{2}{p},\,\,\,b>\frac{1}{2},\,\,\,s>\frac{1}{2}-\frac{2}{p}-\frac{1}{q}\]
the estimate
\begin{equation}\label{s1}
\n{u}{L^p_t(L^q_x)} \leq c \n{u}{\x}
\end{equation}
holds true.
\item[b)] For all $p,\,\,q,\,\,s$ and $b$ satisfying
\[0 \leq \frac{1}{p} \leq \frac{1}{q} \leq \frac{1}{2} \leq \frac{2}{p} + \frac{1}{q} \leq \frac{3}{2},\,\,\,s>0 \,\,\,\mbox{and}\,\,\,b>\frac{3}{4}-\frac{1}{p}-\frac{1}{2q}\]
the estimate (\ref{s1}) is valid.
\item[c)] For all $p,\,\,q,\,\,s$ satisfying
\[0 < \frac{1}{p} \leq \frac{1}{6},\,\,\,0 < \frac{1}{q} \leq \frac{1}{2} - \frac{2}{p},\,\,\,s>\frac{1}{2}-\frac{2}{p}-\frac{1}{q}\]
there exists a $b<\frac{1}{2}$ so that (\ref{s1}) holds true.
\end{itemize}
\end{kor}

Proof: i) By the Sobolev embedding theorem in the time variable we have $\XX{0}{b} \subset L^{\infty}_t (L^2_x)$ for all $b>\frac{1}{2}$. Interpolating this with the above lemma, we obtain (\ref{s1}) whenever $0 \leq \frac{1}{p} \leq \frac{1}{6},\,\,\,s>0$ and $\frac{1}{2} = \frac{2}{p} + \frac{1}{q}$.

ii) Combining this with Sobolev embedding in the space variable, part a) follows. To see part b), one has to interpolate between the result in i) and the trivial case $\XX{0}{0} = L^2_{xt}$.

iii) Now for $p,\,\,q,\,\mbox{and}\,s$ according to the assumptions of part c), there exists a $\theta \in [0,1)$ satisfying
\[\theta \geq 1-\frac{2}{p}\,\,\,\theta>1-\frac{2}{q}\,\,\mbox{and}\,\, s > \frac{3}{2} - \theta -\frac{2}{p}-\frac{1}{q}\,\,.\]
Define $s_1 = \frac{s}{\theta},\,\,\,b_1=\frac{1}{4} +  \frac{1}{4 \theta}$ and $p_1,\,\,q_1$ by $\frac{1}{p}=\frac{1 - \theta}{2} + \frac{\theta}{p_1}$ and $\frac{1}{q}=\frac{1 - \theta}{2} + \frac{\theta}{q_1}$. A simple computation shows, that $p_1,\,\,q_1,\,\,s_1$ and $b_1$ are chosen according to the assumptions of part a). Now part c) with $b= \theta b_1 = \frac{\theta +1}{4} <\frac{1}{2}$ follows by interpolation between this and the trivial case. 
$\hfill \Box$

\vspace{0,5cm}

Next we prove the higherdimensional $L^4$-estimates (cf. \cite{B93}, Prop. 3.6).

\begin{lemma}\label{l23} Let $n \geq 2$. Then for all $s>\frac{n}{2}- \frac{n+2}{4}$ and $b>\frac{1}{2}$ there exists a constant $c = c(s,b)$, so that the following estimate holds:
\[\nf{L^4_{xt}} \leq c \nf{\XX{s}{b}}\,\,.\]
\end{lemma}

Proof: We start by showing that
\[ \sup_{(\xi, \tau) \in \Z^n \times \R} S(\xi, \tau) \leq c N^{4s}\]
for
\begin{eqnarray*}
S(\xi, \tau) & = & \sum_{\xi_1 \in \Z^n}\chi_N (\xi_1)\chi_N (\xi - \xi_1)<\tau + |\xi_1|^2 + |\xi -\xi_1|^2 >^{-2b} \\
& \leq & c \sum_{\xi_1 \in \Z^n}\chi_{2N} (2 \xi_1)\chi_{2N} (2 (\xi - \xi_1))<4 \tau + |2 \xi_1|^2 + |2 (\xi -\xi_1)|^2 >^{-2b}.
\end{eqnarray*}
Here $\chi_N$ denotes the characteristic function of the ball with radius $R$ centered at zero. With $k= 2 \xi_1 - \xi \in \Z^n$ we have
\[k+\xi=2 \xi_1 , \,\,k-\xi=2 (\xi_1 -\xi) \,\,\mbox{and}\,\,|2 \xi_1|^2 + |2 (\xi -\xi_1)|^2 = 2(|\xi| ^2 + |k|^2).\]
Thus we can estimate
\begin{eqnarray*}
S(\xi, \tau) & \leq & c \sum_{k \in \Z^n}\chi_{2N} (k+ \xi)\chi_{2N} (k-\xi)<4 \tau + 2 (|\xi|^2 + |k|^2) >^{-2b} \\
& \leq & c \sum_{k \in \Z^n}\chi_{2N} (k)<2 \tau +  |\xi|^2 + |k|^2 >^{-2b} \\
& = & c \sum_{r \in \N_0}\sum_{k \in \Z^n ,|k|^2 =r} \chi_{4N^2} (r) <2 \tau +  |\xi|^2 + r >^{-2b}\\
& \leq &  c N^{n-2 +2 \epsilon} \sum_{r \in \N_0}<2 \tau +  |\xi|^2 + r >^{-2b} \leq c N^{4s}\,\,,
\end{eqnarray*}
where in the last but one inequality we have used Proposition \ref{p21}. Thus the stated bound on $S(\xi, \tau)$ is proved. Using Lemma 4.2 from \cite{GTV97} again we obtain (with $\int d \nu = \int d \tau_1 \sum_{\xi_1}$ and $(\xi, \tau)=(\xi_1 + \xi_2, \tau_1 + \tau_2)$):
\[\int d \nu \prod_{i=1}^2 \chi_N (\xi_i) <\tau_i + |\xi_i|^2 >^{-2b} \leq c N^{4s}\,\,.\] 
Applying Cauchy-Schwarz and Fubini as in the former proofs we arrive at
\[\n{\int d \nu \prod_{i=1}^2<\tau_i + |\xi_i|^2 >^{-b}f_i(\xi_i,\tau_i)}{L^2_{\xi \tau}} \leq c N^{2s}\prod_{i=1}^2 \n{f_i}{L^2_{\xi \tau}}\]
for all $f_i \in L^2_{\xi \tau}$ which are supported in $\{(\xi, \tau):|\xi|\leq N \}$. Now Lemma \ref{l11} gives for all $u_i \in \XX{0}{b}$, $i=1,2$, having a Fourier transform supported in $\{(\xi, \tau):|\xi|\leq N \}$:
\[\n{u_1 u_2}{L^2_{xt}} \leq c N^{2s} \prod_{i=1}^2 \n{u_i}{ \XX{0}{b}}\,\,.\]
Taking $u=u_1=u_2$ we get
\begin{equation}\label{s3}
\n{u}{L^4_{xt}} \leq c N^s \n{u}{\XX{0}{b}}
\end{equation}
provided the above support condition ist fulfilled.

Now let $(\phi_j)_{j \in \N_0}$ be a smooth partition of the unity. Then, by the Littlewood-Paley-Theorem, we have $\nf{L^4_x(\T^n)} \sim \n{(\sum_{j \in \N_0}|\phi_j * f|^2)^{\frac{1}{2}}}{L^4_x(\T^n)}$. Combining this with the estimate (\ref{s3}) we get
\begin{eqnarray*}
\q{u}{L^4_{xt}} & \leq c & \n{\sum_{j \in \N_0}|\phi_j * u|^2}{L^2_{xt}}\\
& \leq c & \sum_{j \in \N_0} \q{\phi_j * u}{L^4_{xt}} \\
& \leq c & \sum_{j \in \N_0} 2^{2sj}\q{\phi_j * u}{\XX{0}{b}} \leq c \q{u}{\x}
\end{eqnarray*}
$\hfill \Box$

\begin{kor}\label{k22} Let $n \geq 2$:
\begin{itemize}
\item[a)] For all H\"older- and Sobolevexponents $p,\,\,q,\,\,s$ and $b$ satisfying
\[0 \leq \frac{1}{p} \leq \frac{1}{4},\,\,\,0 < \frac{1}{q} \leq \frac{1}{2} - \frac{1}{p},\,\,\,b>\frac{1}{2},\,\,\,s>\frac{n}{2}-\frac{2}{p}-\frac{n}{q}\]
the estimate
\begin{equation}\label{s2}
\n{u}{L^p_t(L^q_x)} \leq c \n{u}{\x}
\end{equation}
holds true.
\item[b)] For all $p,\,\,q,\,\,s$ and $b$ satisfying
\[0 \leq \frac{1}{p} \leq \frac{1}{q} \leq \frac{1}{2} \leq \frac{1}{p} + \frac{1}{q} \leq 1,\,\,\,s> (n-2)(\frac{1}{2}-\frac{1}{q}) \,\,\,\mbox{and}\,\,\,b>1-\frac{1}{p}-\frac{1}{q}\]
the estimate (\ref{s2}) is valid.
\item[c)] For all $p,\,\,q,\,\,s$ satisfying
\[0 < \frac{1}{p} \leq \frac{1}{4},\,\,\,0 < \frac{1}{q} \leq \frac{1}{2} - \frac{1}{p},\,\,\,s>\frac{n}{2}-\frac{2}{p}-\frac{n}{q}\]
there exists a $b<\frac{1}{2}$ so that (\ref{s2}) holds true.
\end{itemize}
\end{kor}

The proof follows the same lines as that of Corollary \ref{k21} and therefore will be omitted.

\vspace{0,5cm}

$Remark:$ Because of $\n{f}{L^p_t(L^q_x)}=\n{\overline{f}}{L^p_t(L^q_x)}$ and $\n{f}{\xm}=\n{\overline{f}}{\x}$ all the results derived in this section so far hold for $\xm$ instead of $\x$. Moreover they are also valid for the corresponding spaces of nonperiodic functions: For $n=1,2$ this is a direct consequence of the Strichartz estimates and \cite{GTV97}, Lemma 2.3. (To obtain Lemma \ref{l21} one has again to interpolate with the trivial case.) For $n \ge 3$, one has to combine Sobolev's embedding theorem with Strichartz and the cited lemma to obtain $\nf{L^4_{xt}} \leq c \nf{\XX{\frac{n-2}{4}}{b}}$.

\begin{lemma}\label{l24} Assume that for some $1< p,q < \infty$ and $s,b \in \R$ the estimate $\n{u}{L^p_t(L^q_x)} \leq c \n{u}{\x}$ is valid. Let $B$ be a ball (or cube) of radius (sidelength) $R$ centered at $\xi_0 \in \Z^n$. Define the projection $P_B u = \Fx^{-1} \chi_B \Fx$, where $\Fx$ is the Fourier transform in the space variable and $\chi_B$ the characteristic function of $B$. Then also the estimate
\[\n{P_B u}{L^p_t(L^q_x)} \leq c R^s\n{u}{\XX{0}{b}}\]
holds true.
\end{lemma}
(cf. \cite{B93}, p.143, (5.6) - (5.8))

Proof: If $\xi_0 =0$, this is obvious. For $\xi_0 \neq 0$ define
\[T_{\xi_0}u(x,t):= \exp{(-ix \xi_0 - it|\xi_0|^2)}u(x+2t\xi_0,t)\,\,.\]
Then $T_{\xi_0}:L^p_t(L^q_x) \r L^p_t(L^q_x)$ is isometric. For the Fourier transform of $T_{\xi_0}u$ the  identity
\[\F T_{\xi_0}u(\xi,\tau) = \F u (\xi+\xi_0,\tau -2 \xi \xi_0 - |\xi_0|^2)\]
is easily checked. Now let $B_0$ be a ball (or cube) of the same size as $B$ centered at zero. Then we have
\begin{eqnarray*}
\F T_{\xi_0} P_B u(\xi,\tau) & = &\F P_B u (\xi+\xi_0,\tau -2 \xi \xi_0 - |\xi_0|^2)\\
& = & \chi_B (\xi+\xi_0)\F u (\xi+\xi_0,\tau -2 \xi \xi_0 - |\xi_0|^2)\\
& = & \chi_{B_0}(\xi)\F T_{\xi_0}u(\xi,\tau) =\F P_{B_0} T_{\xi_0}u(\xi,\tau)\,\,.
\end{eqnarray*}
That is $T_{\xi_0} P_B u = P_{B_0} T_{\xi_0}u$. Moreover, because of
\begin{eqnarray*}
\q{T_{\xi_0}u}{\XX{0}{b}} & = & \int \mu (d \xi) d \tau <\tau+|\xi|^2>^{2b}|\F u (\xi+\xi_0,\tau -2 \xi \xi_0 - |\xi_0|^2)|^2 \\
& = & \int \mu (d \xi) d \tau <\tau+|\xi + \xi_0|^2>^{2b}|\F u (\xi+\xi_0,\tau)|^2 = \q{u}{\XX{0}{b}}
\end{eqnarray*}
$T_{\xi_0}: \XX{0}{b} \r \XX{0}{b}$ is also isometric. Now we can conclude
\begin{eqnarray*}
\n{P_B u}{L^p_t(L^q_x)} & = & \n{T_{\xi_0} P_B u}{L^p_t(L^q_x)} \\
& = & \n{P_{B_0} T_{\xi_0} u}{L^p_t(L^q_x)} \\
& \leq  &  c R^s\n{T_{\xi_0} u}{\XX{0}{b}} = c R^s\n{u}{\XX{0}{b}}
\end{eqnarray*}
$\hfill \Box$

$Remark:$ If $B$ is a ball centered at $\xi_0$ and $-B$ is the ball of the same size centered at $-\xi_0$, then a short computation using $\Fx \overline{u}(\xi) = \overline{\Fx u}(-\xi)$ shows that $P_B \overline{u} = \overline{P_{-B} u}$. From this and $\n{u}{\xm}=\n{\overline{u}}{\x}$ it follows, that Lemma \ref{l24} remains valid with $\x$ replaced by $\xm$. Moreover, as the proof shows, the Lemma is also true in the nonperiodic case.

\section{Multilinear estimates}

\begin{satz}\label{t1} Let $n=1$, $\theta \in (0,\frac{1}{4})$ and $s \geq 0$. Then for all $u_{1,2} \in \XX{s}{\frac{1}{2}}$ supported in $\{(x,t): |t| \leq T\}$ the following estimates are valid:
\begin{itemize}
\item[i)] $\n{\partial_x (\overline{u}_1 \overline{u}_2)}{\XX{s}{-\frac{1}{2}}} \leq c T^{\theta} \n{u_1}{\XX{s}{\frac{1}{2}}}\n{u_2}{\XX{s}{\frac{1}{2}}}$ and
\item[ii)] $\n{\partial_x (\overline{u}_1 \overline{u}_2)}{\y} \leq c T^{\theta} \n{u_1}{\XX{s}{\frac{1}{2}}}\n{u_2}{\XX{s}{\frac{1}{2}}}$
\end{itemize}
\end{satz}

Proof: 1. Preparations: Without loss of generality we may assume $s=0$. Setting $v_i = \overline{u}_i$ the stated inequalities then read
\begin{equation}\label{m1}
\n{\partial_x (v_1 v_2)}{\XX{0}{-\frac{1}{2}}} \leq c T^{\theta} \n{v_1}{\XXm{0}{\frac{1}{2}}}\n{v_2}{\XXm{0}{\frac{1}{2}}}
\end{equation}
and
\begin{equation}\label{m2}
\n{\partial_x (v_1 v_2)}{\yy{0}} \leq c T^{\theta} \n{v_1}{\XXm{0}{\frac{1}{2}}}\n{v_2}{\XXm{0}{\frac{1}{2}}}\,\,.
\end{equation}
To show them, we need the following algebraic inequality:
\begin{eqnarray}\label{ai1}
& & <\xi>^2 + <\xi_1>^2 + <\xi_2>^2 \nonumber \\
& \leq & <\tau + \xi^2> + <\tau_1 - \xi_1^2> + <\tau_2 - \xi_2^2> \\
& \leq & c(<\tau + \xi^2> \chi_A + <\tau_1 - \xi_1^2> + <\tau_2 - \xi_2^2>) \,\,.\nonumber
\end{eqnarray}
Here $A$ denotes the region, where $<\tau + \xi^2> \geq \max_{i=1}^2 <\tau_i - \xi_i^2>$. (For the variables $(\xi,\xi_1,\xi_2)$ and $(\tau,\tau_1,\tau_2)$ we will have $\xi=\xi_1+\xi_2$ and $\tau=\tau_1+\tau_2$ throughout this proof.) Defining $f_i(\xi, \tau)= <\tau - \xi^2>^{\frac{1}{2}} \F v_i(\xi, \tau)$ for $i=1,2$ we have \\ $\n{v_i}{\XXm{0}{\frac{1}{2}}}=\n{f_i}{L^2_{\xi,\tau}}$. Now, for given $\theta \in (0,\frac{1}{4})$ we fix $\epsilon = \frac{1}{4} (\frac{1}{4}-\theta)$.

2. Estimation of (\ref{m1}): By Lemma \ref{l11} and (\ref{ai1}) we have:
\begin{eqnarray*}
& & \n{\partial_x (v_1 v_2)}{\XXm{0}{-\frac{1}{2}}} \\
&=& c \n{<\tau + \xi^2>^{-\frac{1}{2}}\xi \int \mu (d\xi_1) d \tau_1 \prod_{i=1}^2 <\tau_i - \xi_i^2>^{-\frac{1}{2}}f_i(\xi_i,\tau_i)}{L^2_{\xi,\tau}} \\
& \leq & c \sum_{i=1}^3 N_i
\end{eqnarray*}
with
\[N_1 = \n{ \int \mu (d\xi_1) d \tau_1 \prod_{i=1}^2 <\tau_i - \xi_i^2>^{-\frac{1}{2}}f_i(\xi_i,\tau_i)}{L^2_{\xi,\tau}},\]
\[N_2= \n{<\tau + \xi^2>^{-\frac{1}{2}} \int \mu (d\xi_1) d \tau_1 <\tau_2 - \xi_2^2>^{-\frac{1}{2}}\prod_{i=1}^2 f_i(\xi_i,\tau_i)}{L^2_{\xi,\tau}}\]
and
\[N_3= \n{<\tau + \xi^2>^{-\frac{1}{2}} \int \mu (d\xi_1) d \tau_1 <\tau_1 - \xi_1^2>^{-\frac{1}{2}}\prod_{i=1}^2 f_i(\xi_i,\tau_i)}{L^2_{\xi,\tau}}\,\,.\]
Lemma \ref{l11}, H\"olders inequality, Lemma \ref{l21} and Lemma \ref{l12} are now applied to obtain
\begin{eqnarray*}
N_1 = \n{v_1 v_2}{L^2_{x,t}} & \leq & \n{v_1}{L^4_{x,t}}\n{v_2}{L^4_{x,t}} \\
& \leq & c \n{v_1}{\XXm{0}{\frac{3}{8} + \epsilon}}\n{v_2}{\XXm{0}{\frac{3}{8} + \epsilon}} \\
& = & c \n{\psi _{2T}v_1}{\XXm{0}{\frac{3}{8} + \epsilon}}\n{\psi _{2T}v_2}{\XXm{0}{\frac{3}{8} + \epsilon}} \\
& \leq & c T^{\frac{1}{4}- 4\epsilon} \n{v_1}{\XXm{0}{\frac{1}{2} - \epsilon}}\n{v_2}{\XXm{0}{\frac{1}{2} - \epsilon}}\,\,.
\end{eqnarray*}
Similarly we get
\begin{eqnarray*}
N_2 = \n{(\F ^{-1}f_1) v_2}{\XX{0}{-\frac{1}{2}}} & \leq & \n{ \psi _{2T} (\F ^{-1}f_1) v_2}{\XX{0}{-\frac{1}{2} + \epsilon}} \\
& \leq & c T^{\frac{1}{8}- 2\epsilon}\n{(\F ^{-1}f_1) v_2}{\XX{0}{-\frac{3}{8} - \epsilon}}\\
& \leq & c T^{\frac{1}{8}- 2\epsilon}\n{(\F ^{-1}f_1) v_2}{L^{\frac{4}{3}}_{x,t}} \\
& \leq & c T^{\frac{1}{8}- 2\epsilon}\n{\F ^{-1}f_1}{L^2_{x,t}}\n{v_2}{L^4_{x,t}}\\
& \leq & c T^{\frac{1}{8}- 2\epsilon} \n{v_1}{\XXm{0}{\frac{1}{2}}}\n{\psi _{2T}v_2}{\XXm{0}{\frac{3}{8} + \epsilon}} \\
& \leq & c T^{\frac{1}{4}- 4\epsilon} \n{v_1}{\XXm{0}{\frac{1}{2}}}\n{v_2}{\XXm{0}{\frac{1}{2}}} \,\,.
\end{eqnarray*}
By exchanging $v_1$ and $v_2$ we get the same upper bound for $N_3$. So, because of $\theta = \frac{1}{4}- 4\epsilon$, the estimate (\ref{m1}) is proved.

3. Estimation of (\ref{m2}): Using Lemma \ref{l11} and (\ref{ai1}) we get
\begin{eqnarray*}
& & \n{\partial_x (v_1 v_2)}{\yy{0}} \\
&=& c \n{<\tau + \xi^2>^{-1}\xi \int \mu (d\xi_1) d \tau_1 \chi_A \prod_{i=1}^2 <\tau_i - \xi_i^2>^{-\frac{1}{2}}f_i(\xi_i,\tau_i)}{L^2_{\xi}(L^1_{\tau})} \\
& \leq & c \sum_{i=1}^3 N_i\,\,\,,
\end{eqnarray*}
where
\[N_1 = \n{<\tau + \xi^2>^{-\frac{1}{2}} \int \mu (d\xi_1) d \tau_1 \chi_A \prod_{i=1}^2 <\tau_i - \xi_i^2>^{-\frac{1}{2}}f_i(\xi_i,\tau_i)}{L^2_{\xi}(L^1_{\tau})}\,\,,\]
\[N_2= \n{<\tau + \xi^2>^{-1} \int \mu (d\xi_1) d \tau_1 <\tau_2 - \xi_2^2>^{-\frac{1}{2}}\prod_{i=1}^2 f_i(\xi_i,\tau_i)}{L^2_{\xi}(L^1_{\tau})}\]
and
\[N_3= \n{<\tau + \xi^2>^{-1} \int \mu (d\xi_1) d \tau_1 <\tau_1 - \xi_1^2>^{-\frac{1}{2}}\prod_{i=1}^2 f_i(\xi_i,\tau_i)}{L^2_{\xi}(L^1_{\tau})}\,\,.\]
In order to estimate $N_1$ we define
\[g_i(\xi,\tau):=<\tau - \xi^2>^{\frac{3}{8} + \epsilon} \F v_i(\xi, \tau)=<\tau - \xi^2>^{-\frac{1}{8} + \epsilon}f_i(\xi,\tau)\,\,.\]
Then it is $\n{g_i}{L^2_{\xi,\tau}}=\n{v_i}{\XXm{0}{\frac{3}{8} + \epsilon}}$ and
\[N_1 = \n{<\tau + \xi^2>^{-\frac{1}{2}} \int \mu (d\xi_1) d \tau_1 \chi_A \prod_{i=1}^2 <\tau_i - \xi_i^2>^{-\frac{3}{8} - \epsilon}g_i(\xi_i,\tau_i)}{L^2_{\xi}(L^1_{\tau})}\,\,.\]
Since in $A$ we have $<\tau + \xi^2> \geq \max_{i=1}^2 <\tau_i - \xi_i^2>$ as well as $<\tau + \xi^2> \geq c< \xi_1>^2$, we obtain
\[N_1 \leq c \n{ \int \mu (d\xi_1) d \tau_1 <\xi_1>^{-\frac{1}{2} - 2 \epsilon} \prod_{i=1}^2 <\tau_i - \xi_i^2>^{-\frac{1+\epsilon}{2}}g_i(\xi_i,\tau_i)}{L^2_{\xi}(L^1_{\tau})}\,\,,\]
which we shall now estimate by duality. Therefor let $f_0 \in L^2_{\xi}$ with
$\n{f_0}{L^2_{\xi}}=1$ and $f_0 \geq 0$. Now applying Cauchy-Schwarz' inequality first in the $\tau$- and then in the $\xi$-variables we get the desired upper bound for $N_1$:
\begin{eqnarray*}
& & \int \mu (d\xi d\xi_1) d \tau d \tau_1 f_0(\xi)<\xi_1>^{-\frac{1}{2} - 2 \epsilon} \prod_{i=1}^2 <\tau_i - \xi_i^2>^{-\frac{1+\epsilon}{2}}g_i(\xi_i,\tau_i) \\
&=&\int \mu (d\xi_1 d\xi_2) d \tau_1 d \tau_2 f_0(\xi_1 + \xi_2)<\xi_1>^{-\frac{1}{2} - 2 \epsilon} \prod_{i=1}^2 <\tau_i - \xi_i^2>^{-\frac{1+\epsilon}{2}}g_i(\xi_i,\tau_i) \\
&\leq c & \int \mu (d\xi_1 d\xi_2)f_0(\xi_1 + \xi_2)<\xi_1>^{-\frac{1}{2} - 2 \epsilon}\prod_{i=1}^2 (\int d \tau_i|g_i(\xi_i,\tau_i)|^2)^{\frac{1}{2}} \\
&\leq c &\prod_{i=1}^2\n{g_i}{L^2_{\xi,\tau}}
\leq c \prod_{i=1}^2\n{v_i}{\XXm{0}{\frac{3}{8} + \epsilon}} 
\leq c T^{\frac{1}{4}- 
4\epsilon}\prod_{i=1}^2\n{v_i}{\XXm{0}{\frac{1}{2}}}\,\,,
\end{eqnarray*}
where in the last step we have used Lemma \ref{l12} from the introduction. To estimate $N_2$ we apply Cauchy-Schwarz on $\int d \tau$:
\begin{eqnarray*}
N_2 \leq c \n{<\tau + \xi^2>^{-\frac{1}{2} +\epsilon} \int \mu (d\xi_1) d \tau_1 <\tau_2 - \xi_2^2>^{-\frac{1}{2}}\prod_{i=1}^2 f_i(\xi_i,\tau_i)}{L^2_{\xi,\tau}} \\ =  
\n{ \psi _{2T} (\F ^{-1}f_1) v_2}{\XX{0}{-\frac{1}{2} + \epsilon}} \,\,.\hspace{4cm}
\end{eqnarray*}
This was already shown to be bounded by
\[c T^{\frac{1}{4}- 4\epsilon}\prod_{i=1}^2\n{v_i}{\XXm{0}{\frac{1}{2}}}\,\,.\]
The same upper bound for $N_3$ is obtained by exchanging $v_1$ and $v_2$, so the estimate (\ref{m2}) is proved, too.
$\hfill \Box$

\begin{satz}\label{t2} Let $n, m \in \N$ with $m \ge 2$ and $m+n \ge 4$. Assume in addition, that $s > \frac{n}{2} -\frac{1}{m-1}$. Then there exists a $\theta > 0$, so that for all $0<T \le 1$ and for all $u_i \in \XX{s}{\frac{1}{2}},\,\,\,1\le i \le m$ having support in $\{(x,t): |t|<T\}$ the estimates
\begin{itemize}
\item[i)] $\n{ \prod_{i=1}^m\overline{u}_i}{\XX{s+1}{-\frac{1}{2}}} \leq c T^{\theta} \prod_{i=1}^m \n{u_i}{\XX{s}{\frac{1}{2}}}$ and
\item[ii)] $\n{ \prod_{i=1}^m\overline{u}_i}{\yy{s+1}} \leq c T^{\theta} \prod_{i=1}^m \n{u_i}{\XX{s}{\frac{1}{2}}}$
\end{itemize}
hold.
\end{satz}

Before we prove the theorem, we must introduce some notation and derive some preparatory lemmas. First, for a subset $M \subset \R^n$ or $M \subset \Z^n$, we define the projections $P_M := \Fx ^{-1} \chi_M \Fx$, where $\chi_M$ denotes the characteristic function of the set $M$. Especially we require for $l \in \N_0$:
\begin{itemize}
\item $P_l := P_{B_{2^l}}$ for the (closed) ball $B_{2^l}$ of radius $2^l$ centered at zero ($P_{-1} = 0$),
\item $P_{\Delta l} := P_l - P_{l-1}$, as well as
\item $P_{Q^l_{\alpha}}$, where $\alpha \in \Z^n$ and $Q^l_{\alpha}$ is a cube  of sidelength $2^l$ centered at $2^l \alpha$, so that
\[\R^n = \sum_{\alpha \in \Z^n}Q^l_{\alpha} \hspace{1cm} \mbox{respectively} \hspace{1cm} \Z^n = \sum_{\alpha \in \Z^n}Q^l_{\alpha}\,\,.\]
\end{itemize}

Next we shall fix a couple of H\"older- and Sobolevexponents to be used below:

1. We choose $\frac{1}{p}=\frac{1}{(n+2)(m-1)}$. Then for any $s > \frac{n}{2} -\frac{1}{m-1}$ by corollaries \ref{k21} and \ref{k22}, part c), there exists a $b<\frac{1}{2}$, so that the following estimate holds:
\begin{equation}\label{m3}
\n{u}{L^p_{xt}} \leq c \n{u}{\xpm}
\end{equation}

2. Next we have $\frac{1}{p_0}=\frac{1}{6} + \epsilon$ for $n=1$ respectively $\frac{1}{p_0}=\frac{1}{4} + \epsilon$ for $n \ge 2$ and $s_0 = \epsilon$ if $n=1$ respectively $s_0 = (n-2)(\frac{1}{2}-\frac{1}{p_0}) +\epsilon = \frac{n-2}{4} + (3-n)\epsilon$ if $n \ge 2$. Then, if $\epsilon > 0$ is chosen appropriately small, by corollaries \ref{k21} and \ref{k22}, part b), and Lemma \ref{l24} there exists a $b < \frac{1}{2}$ for which we have the estimate
\begin{equation}\label{m4}
\n{P_B u}{L^{p_0}_{xt}} \leq c R^{s_0}\n{u}{\XXpm{0}{b}}\,\,,
\end{equation}
whenever $B$ is a ball or cube of size $R$. Dualizing the last inequality, we obtain
\begin{equation}\label{m5}
\n{P_B u}{\XXpm{0}{-b}}\leq c R^{s_0}\n{u}{L^{p'_0}_{xt}}\,\,,
\end{equation}
where $\frac{1}{p'_0}=\frac{5}{6} - \epsilon$ for $n=1$ respectively $\frac{1}{p'_0}=\frac{3}{4} - \epsilon$ for $n \ge 2$.

3. We choose $\frac{1}{p_1}=\frac{1}{3} - \epsilon -\frac{m-2}{3(m-1)}$ for $n=1$ respectively $\frac{1}{p_1}=\frac{1}{4} - \epsilon -\frac{m-2}{(n+2)(m-1)}$ for $n \ge 2$ and $s_1 = \frac{n}{2}-\frac{n+2}{p_1}+\epsilon$. Then it is $s_1 = \frac{1}{2}-\frac{1}{m-1}+ 4\epsilon$  if $n=1$ respectively $s_1 =\frac{n+2}{4} -\frac{1}{m-1} + (n+3)\epsilon$ if $n \ge 2$, and by corollaries \ref{k21}, \ref{k22}, part c), and Lemma \ref{l24} there exists a $b < \frac{1}{2}$ for which 
\begin{equation}\label{m6}
\n{P_B u}{L^{p_1}_{xt}} \leq c R^{s_1}\n{u}{\XXpm{0}{b}}\,\,.
\end{equation}
Observe that our choice guarantees
\[\frac{1}{p_0} + \frac{1}{p_1} + \frac{m-2}{p} = \frac{1}{2} \hspace{1cm}\mbox{resp.}\hspace{1cm} \frac{1}{p_1} + \frac{1}{2} + \frac{m-2}{p} = \frac{1}{p'_0}\]
for the H\"older applications as well as for $\epsilon$ sufficiently small
$s_0 + s_1 -s < 0$.

For $m \ge 3$ in addition we shall need the following parameters:

4. Assuming $\frac{s}{n}<\frac{1}{2}$ without loss of generality, we may choose $\frac{1}{q}=\frac{1}{2}-\frac{s}{n}>0$, so that the Sobolev embedding $H_x^s \subset L^q_x$ holds.

5. In the case of space dimension $n=1$ we define $\frac{1}{r_0}=\frac{1}{6} - \frac{m-3}{6(m-1)} - \epsilon$, $\frac{1}{q_0}= s + \frac{1}{6} - \frac{2(m-3)}{3(m-1)} - \epsilon$ and $\sigma_1 =\epsilon$, if $m=3$, as well as $\sigma_1 =\frac{1}{2}-\frac{2}{r_0}-\frac{1}{q_0}+\epsilon = \frac{m-3}{m-1} - s + 4\epsilon$ if $m \ge 4$. For $n \ge 2$ let $\frac{1}{r_0}=\frac{1}{4} - \frac{m-3}{(n+2)(m-1)} - 2\epsilon$, $\frac{1}{q_0}= \frac{s}{n} - \frac{1}{4} - \frac{m-3}{(n+2)(m-1)} - \epsilon$ and $\sigma_1 = \frac{n}{2} - \frac{2}{r_0} - \frac{n}{q_0} +\epsilon = \frac{3n}{4} + \frac{1}{2} -\frac{2}{m-1} - s + (n+5) \epsilon$. Then, for some $b<\frac{1}{2}$, we have the estimate
\begin{equation}\label{m7}
\n{P_B u}{L^{r_0}_t(L^{q_0}_x)} \leq c R^{\sigma_1}\n{u}{\XXpm{0}{b}}\,\,.
\end{equation}
In general, this follows from part c) of the corollaries \ref{k21}, \ref{k22}, except in the case $n=1,\,\,m=3$, where one can use part b) of corollary \ref{k21}. (Here we assume $s \le \frac{1}{3}$ in the cases $n=1\,,\,m \in \{3,4\}$.)

6. We close our list of parameters by choosing $\frac{1}{r_1}=\frac{1}{6} - \frac{m-3}{6(m-1)}$, $\frac{1}{q_1} = \frac{1}{2}-\frac{2}{r_1} = \frac{1}{6} + \frac{m-3}{3(m-1)}$ for $n=1$ respectively $\frac{1}{r_1}=\epsilon$, $\frac{1}{q_1} = \frac{1}{2}$ for $n \ge 2$. Then, by corollary \ref{k21}, part c), in the case of space dimension $n=1$ and by Sobolev embedding in the time variable in the case of $n \ge 2$, we have the estimate
\begin{equation}\label{m8}
\n{P_B u}{L^{r_1}_t(L^{q_1}_x)} \leq c R^{\epsilon}\n{u}{\XXpm{0}{b}}
\end{equation}
for some $b<\frac{1}{2}$. Now for the H\"older applications we have
\[\frac{1}{r_0} + \frac{1}{2} + \frac{1}{r_1} + \frac{m-3}{p} =  \frac{1}{q_0} + \frac{1}{q} + \frac{1}{q_1} + \frac{m-3}{p} = \frac{1}{p'_0}\]
as well as for $\epsilon$ sufficiently small
$s_0 + \sigma_1 + \epsilon -s < 0$.

\begin{lemma}\label{l31} Let $n,m \in \N$ with $m \ge 2$ and $n+m \ge 4$. Then for $s > \frac{n}{2} -\frac{1}{m-1}$ there exists a $b<\frac{1}{2}$, so that for all $v_i \in \xm$, $1 \leq i,j \leq m$ the following estimate is valid:
\[\n{(J^s v_j)\prod^m_{i=1, i \neq j} v_i}{L^2_{xt}} \leq c \prod_{i=1}^m \n{v_i}{\xm}\,\,,\]
where $J^s = \Fx ^{-1} <\xi >^s \Fx$.
\end{lemma}

Proof: Writing
\[\prod^m_{i=1 \atop i \neq j} v_i = \lim_{n \in \N_0} \prod^m_{i=1 \atop i \neq j} P_l v_i =\sum_{l \in \N_0} (\prod^m_{i=1 \atop i \neq j} P_l v_i - \prod^m_{i=1 \atop i \neq j} P_{l-1} v_i)\,\,,\]
where
\[\prod^m_{i=1 \atop i \neq j} P_l v_i - \prod^m_{i=1 \atop i \neq j} P_{l-1} v_i = \sum_{k=1 \atop k \neq j}^m (\prod_{i<k \atop i \neq j} P_{l-1} v_i) P_{\Delta l}v_k (\prod_{i>k \atop i \neq j} P_l v_i)\,\,,\]
we obtain
\begin{eqnarray}\label{m9}
& & \n{(J^s v_j)\prod_{ i \neq j} v_i}{L^2_{xt}}  \nonumber \\
& \leq & \sum_{l \in \N_0} \sum_{k=1 \atop k \neq j}^m \n{(J^s v_j)(\prod_{i<k \atop i \neq j} P_{l-1}v_i P_{\Delta l}v_k (\prod_{i>k \atop i \neq j} P_l v_i))}{L^2_{xt}} \\
& \leq & \sum_{l \in \N_0} \sum_{k=1 \atop k \neq j}^m \n{(J^s v_j) (P_{\Delta l}v_k)      (\prod_{ i \neq j} P_l v_i)}{L^2_{xt}} \,\,.\nonumber
\end{eqnarray}
Next we estimate the contribution for fixed $l$ and $k$:
\begin{eqnarray*}
& & \q{(J^s v_j) (P_{\Delta l}v_k) (\prod_{ i \neq j} P_l v_i)}{L^2_{xt}} \\
& = & \q{\sum_{\alpha \in \Z^n}(P_{Q^l_{\alpha}} J^s v_j) (P_{\Delta l}v_k) (\prod_{ i \neq j} P_l v_i)}{L^2_{xt}} \\
& = & \sum_{\alpha , \beta \in \Z^n} <(P_{Q^l_{\alpha}} J^s v_j) (P_{\Delta l}v_k) (\prod_{ i \neq j} P_l v_i),(P_{Q^l_{\beta}} J^s v_j) (P_{\Delta l}v_k) (\prod_{ i \neq j} P_l v_i)>
\end{eqnarray*}

Now the sequence $\{(P_{Q^l_{\alpha}} J^s v_j) (P_{\Delta l}v_k) (\prod_{ i \neq j} P_l v_i)\}_{\alpha \in \Z^n}$ is almost orthogonal in the following sense: The support of $\F (P_{\Delta l}v_k) (\prod_{ i \neq j} P_l v_i)$ is contained in $\{(\xi,\tau): |\xi| \le (m-1)2^l\}$, and thus $\F (P_{Q^l_{\alpha}} J^s v_j)(P_{\Delta l}v_k) (\prod_{ i \neq j} P_l v_i)$ is supported in $C \times  \R$, where $C$ is a cube centered at $2^l \alpha$ having the sidelength $m 2^l$. So for $|2^l \alpha - 2^l \beta|>c_n 2^l m$, that is for $|\alpha - \beta|>c_n m$, the above expressions are disjointly supported. Thus for these values of $\alpha$ and $\beta$ we do not get any contribution to the last sum, which we now can estimate by
\begin{eqnarray}\label{m10}
\sum_{\alpha \in \Z^n}\sum_{\beta \in \Z^n \atop |\beta| \le c_n m} <(P_{Q^l_{\alpha}} J^s v_j) (P_{\Delta l}v_k) (\prod_{ i \neq j} P_l v_i),(P_{Q^l_{\alpha + \beta}} J^s v_j) (P_{\Delta l}v_k) (\prod_{ i \neq j} P_l v_i)> \nonumber \\ 
 \le  c \sum_{\alpha \in \Z^n} \q{(P_{Q^l_{\alpha}} J^s v_j) (P_{\Delta l}v_k) (\prod_{ i \neq j} P_l v_i)}{L^2_{xt}} \hspace{3cm}\\
\le c \sum_{\alpha \in \Z^n} \q{(P_{Q^l_{\alpha}} J^s v_j) (P_{\Delta l}v_k) (\prod_{ i \neq j} v_i)}{L^2_{xt}}\,\,.\hspace{3.2cm} \nonumber
\end{eqnarray}

Next we use H\"older's inequality, (\ref{m3}), (\ref{m4}) and (\ref{m6}) to get
\begin{eqnarray}\label{m11}
& & \n{(P_{Q^l_{\alpha}} J^s v_j) (P_{\Delta l}v_k) (\prod_{ i \neq j} v_i)}{L^2_{xt}} \nonumber \\
& \le & \n{P_{Q^l_{\alpha}} J^s v_j}{L^{p_0}_{xt}}\n{P_{\Delta l}v_k}{L^{p_1}_{xt}} \prod_{i \neq k,j}\n{v_i}{L^p_{xt}} \\
& \le & c 2^{l(s_0 + s_1)} \n{P_{Q^l_{\alpha}} J^s v_j}{\XXm{0}{b}}\n{P_{\Delta l}v_k}{\XXm{0}{b}} \prod_{i \neq k,j}\n{v_i}{\xm} \nonumber
\end{eqnarray}
for some $b<\frac{1}{2}$. Using $\n{P_{\Delta l}v_k}{\XXm{0}{b}} \le c 2^{-sl}\n{v_k}{\XXm{s}{b}}$ we combine (\ref{m10}) and (\ref{m11}) to obtain:
\begin{eqnarray*}
& & \q{(J^s v_j) (P_{\Delta l}v_k) (\prod_{ i \neq j} P_l v_i)}{L^2_{xt}} \\
& \le & c 2^{2l(s_0 + s_1 - s)} \sum_{\alpha \in \Z^n} \q{P_{Q^l_{\alpha}} J^s v_j}{\XXm{0}{b}} \prod_{i \neq j}\q{v_i}{\xm} \\
& = & c 2^{2l(s_0 + s_1 - s)} \prod^m_{i=1}\q{v_i}{\xm}\,\,.
\end{eqnarray*}
Inserting the square root of this into (\ref{m9}) and summing up over $k$ and $l$ we can finish the proof.
$\hfill \Box$

\begin{kor}\label{k31} For $n,m$ and $s$ as in the previous lemma there exists a $b< \frac{1}{2}$, so that for all $v_i \in \XXm{s}{\frac{1}{2}}$, $1 \le i,j \le m$ the following estimate holds true:
\[\n{(\Lambda ^{\frac{1}{2}}J^s v_j)\prod^m_{i=1, i \neq j} v_i}{\XX{0}{-b}} \leq c \n{v_j}{\XXm{s}{\frac{1}{2}}} \prod_{i=1 \atop i \neq j}^m \n{v_i}{\xm}\,\,,\]
where $\Lambda ^{\frac{1}{2}} = \F^{-1}<\tau -|\xi|^2>^{\frac{1}{2}}\F$.
\end{kor}

Proof: Let the $v_i$'s be fixed for $i \neq j$. Then the previous lemma tells us, that the linear mapping
\[A_j : \xm \r L^2_{xt}, \hspace{1cm}f \mapsto (J^s f)\prod^m_{i=1 \atop i \neq j} v_i\]
is bounded with norm $\|A_j\| \le c \prod_{i=1 \atop i \neq j}^m \n{v_i}{\xm}$. The adjoint mapping $A_j^*$, given by
\[A_j^* :  L^2_{xt} \r \XXm{-s}{-b}, \hspace{1cm}g \mapsto J^s (g\prod^m_{i=1 \atop i \neq j} \overline{v_i})\]
then is also bounded with $\|A_j^*\|=\|A_j\|$. From this we get for $g=\overline{\Lambda ^{\frac{1}{2}}J^s v_j}$:
\begin{eqnarray*}
 \n{(\Lambda ^{\frac{1}{2}}J^s v_j)\prod^m_{i=1, i \neq j} v_i}{\XX{0}{-b}} 
& = & \n{J^s(\overline{\Lambda ^{\frac{1}{2}}J^s v_j})\prod^m_{i=1, i \neq j} \overline{v_i}}{\XXm{-s}{-b}} \\
\le  c \n{\Lambda ^{\frac{1}{2}}J^s v_j}{L^2_{xt}}\prod_{i=1 \atop i \neq j}^m \n{v_i}{\xm} 
 & = & c \n{v_j}{\XXm{s}{\frac{1}{2}}}\prod_{i=1 \atop i \neq j}^m \n{v_i}{\xm}
\end{eqnarray*}
$\hfill \Box$

\begin{lemma}\label{l32} Let $n,m \in \N$ with $m \ge 2$, $n+m \ge 4$ and $s \in (\frac{n}{2} -\frac{1}{m-1},\frac{n}{2})$. For $n=1$, $m \in \{3,4\}$ assume in addition, that $s \le \frac{1}{3}$. Then there exists a $b<\frac{1}{2}$, so that for all $v_i \in \XXm{s}{\frac{1}{2}}$, $1 \leq i,j \leq m$ the following estimate is valid:
\[\n{(J^s v_i)(\Lambda ^{\frac{1}{2}}v_j)\prod^m_{k=1, k \neq i,j} v_k}{\XX{0}{-b}} \leq c \n{v_j}{\XXm{s}{\frac{1}{2}}}\prod_{k=1 \atop k \neq j}^m \n{v_k}{\xm}\]
Here again we have $\Lambda ^{\frac{1}{2}} = \F^{-1}<\tau -|\xi|^2>^{\frac{1}{2}}\F$.
\end{lemma}

Proof: 1. Similarly as in the proof of the previous lemma we write
\[\Lambda ^{\frac{1}{2}}v_j \prod^m_{k=1 \atop k \neq i,j} v_k = \sum_{l \in \N_0 } ( P_l \Lambda ^{\frac{1}{2}}v_j \prod^m_{k=1 \atop k \neq i,j} P_l v_k - P_{l-1}\Lambda ^{\frac{1}{2}}v_j \prod^m_{k=1 \atop k \neq i,j} P_{l-1} v_k )\]
with
\begin{eqnarray*}
& & P_l \Lambda ^{\frac{1}{2}}v_j \prod^m_{k=1 \atop k \neq i,j} P_l v_k - P_{l-1}\Lambda ^{\frac{1}{2}}v_j \prod^m_{k=1 \atop k \neq i,j} P_{l-1} v_k \\
& = & P_{\Delta l}\Lambda ^{\frac{1}{2}}v_j \prod^m_{k=1 \atop k \neq i,j} P_l v_k + P_{l-1}\Lambda ^{\frac{1}{2}}v_j \sum_{k \neq i,j}(\prod_{\nu <k \atop \nu \neq i,j} P_{l-1} v_{\nu}) P_{\Delta l}v_k (\prod_{\nu > k \atop \nu \neq i,j} P_l v_{\nu})\,\,.
\end{eqnarray*}
From this we obtain for arbitrary $b$:
\begin{eqnarray}\label{m12}
& & \n{(J^s v_i)(\Lambda ^{\frac{1}{2}}v_j)\prod_{ k \neq i,j} v_k}{\XX{0}{-b}} \nonumber \\
& \le & \sum_{l \in \N_0} \n{(J^s v_i)(P_{\Delta l}\Lambda ^{\frac{1}{2}}v_j)\prod_{ k \neq i,j} P_l v_k}{\XX{0}{-b}} \\
& + & \sum_{k \neq i,j} \sum_{l \in \N_0} \n{(J^s v_i)(P_l \Lambda ^{\frac{1}{2}}v_j) (P_{\Delta l} v_k) \prod_{\nu  \neq i,j,k} P_l v_{\nu}}{\XX{0}{-b}} \nonumber
\end{eqnarray}

2. Next we show that for some $b<\frac{1}{2}$ the estimate
\begin{equation}\label{m13}
\n{(J^s v_i)(P_{\Delta l}\Lambda ^{\frac{1}{2}}v_j)\prod_{ k \neq i,j} P_l v_k}{\XX{0}{-b}} \le c 2^{l(s_0 + s_1 -s)} \n{v_j}{\XXm{s}{\frac{1}{2}}}\prod_{i=1 \atop i \neq j}^m \n{v_i}{\xm}
\end{equation}
holds true. To see this, we start from
\begin{eqnarray*}
& & \q{(J^s v_i)(P_{\Delta l}\Lambda ^{\frac{1}{2}}v_j)\prod_{ k \neq i,j} P_l v_k}{\XX{0}{-b}} \\
& = & \q{\sum_{\alpha \in \Z^n}(P_{Q^l_{\alpha}} J^s v_i)(P_{\Delta l}\Lambda ^{\frac{1}{2}}v_j)\prod_{ k \neq i,j} P_l v_k}{\XX{0}{-b}} \\
& \le & c \sum_{\alpha \in \Z^n}\q{(P_{Q^l_{\alpha}} J^s v_i)(P_{\Delta l}\Lambda ^{\frac{1}{2}}v_j)\prod_{ k \neq i,j} P_l v_k}{\XX{0}{-b}}\,\,,
\end{eqnarray*}
where in the last step we have used the almost orthogonality of the sequence $\{(P_{Q^l_{\alpha}} J^s v_i)(P_{\Delta l}\Lambda ^{\frac{1}{2}}v_j)\prod_{ k \neq i,j} P_l v_k\}_{\alpha \in \Z^n}$. Now we use (\ref{m5}), H\"olders inequality, (\ref{m6}) and (\ref{m3}) to obtain for some $b<\frac{1}{2}$
\begin{eqnarray*}
& &\n{(P_{Q^l_{\alpha}} J^s v_i)(P_{\Delta l}\Lambda ^{\frac{1}{2}}v_j)\prod_{ k \neq i,j} P_l v_k}{\XX{0}{-b}} \\
& \le & c 2^{l s_0} \n{(P_{Q^l_{\alpha}} J^s v_i)(P_{\Delta l}\Lambda ^{\frac{1}{2}}v_j)\prod_{ k \neq i,j} P_l v_k}{L^{p'_0}_{xt}} \\
& \le & c 2^{l s_0} \n{P_{Q^l_{\alpha}} J^s v_i}{L^{p_1}_{xt}} \n{P_{\Delta l}\Lambda ^{\frac{1}{2}}v_j}{L^2_{xt}} \prod_{ k \neq i,j}\n{P_l v_k}{L^p_{xt}} \\
& \le & c 2^{l (s_0 + s_1)}\n{P_{Q^l_{\alpha}} J^s v_i}{\XXm{0}{b}} \n{P_{\Delta l}\Lambda ^{\frac{1}{2}}v_j}{L^2_{xt}} \prod_{ k \neq i,j}\n{ v_k}{\xm}\,\,.
\end{eqnarray*}
Using $\n{P_{\Delta l}\Lambda ^{\frac{1}{2}}v_j}{L^2_{xt}} \le c 2^{-ls} \n{v_j}{\XXm{s}{\frac{1}{2}}}$ we get
\begin{eqnarray*}
\q{(P_{Q^l_{\alpha}} J^s v_i)(P_{\Delta l}\Lambda ^{\frac{1}{2}}v_j)\prod_{ k \neq i,j} P_l v_k}{\XX{0}{-b}} \\
\le c 2^{2l(s_0 + s_1 -s)}\q{P_{Q^l_{\alpha}} J^s v_i}{\XXm{0}{b}}\q{v_j}{\XXm{s}{\frac{1}{2}}}\prod_{ k \neq i,j}\q{ v_k}{\xm} \,\,.
\end{eqnarray*}
Now summing up over $\alpha$ we arrive at the square of (\ref{m13}).

3. Now we show that there exists a $b<\frac{1}{2}$ for which
\begin{equation}\label{m14}
\n{(J^s v_i)(P_l \Lambda ^{\frac{1}{2}}v_j) (P_{\Delta l} v_k) \prod_{\nu  \neq i,j,k} P_l v_{\nu}}{\XX{0}{-b}} \le c 2^{l(s_0 + \sigma_1 + \epsilon -s)} \n{v_j}{\XXm{s}{\frac{1}{2}}}\prod_{i=1 \atop i \neq j}^m \n{v_i}{\xm}\,\,.
\end{equation}
Therefor again we write $J^s v_i = \sum_{\alpha \in \Z^n} P_{Q^l_{\alpha}} J^s v_i$ and use the almost orthogonality of $\{(P_{Q^l_{\alpha}} J^s v_i)(P_l \Lambda ^{\frac{1}{2}}v_j) (P_{\Delta l} v_k) \prod_{\nu  \neq i,j,k} P_l v_{\nu}\}_{\alpha \in \Z^n}$ to obtain
\begin{eqnarray*}
& & \q{(J^s v_i)(P_l \Lambda ^{\frac{1}{2}}v_j) (P_{\Delta l} v_k) \prod_{\nu  \neq i,j,k} P_l v_{\nu}}{\XX{0}{-b}} \\
& \le & c \sum_{\alpha \in \Z^n}\q{(P_{Q^l_{\alpha}}J^s v_i)(P_l \Lambda ^{\frac{1}{2}}v_j) (P_{\Delta l} v_k) \prod_{\nu  \neq i,j,k} P_l v_{\nu}}{\XX{0}{-b}}\,\,.
\end{eqnarray*}
Then we use (\ref{m5}), H\"olders inequality, (\ref{m7}), Sobolev embedding in $x$, (\ref{m8}) and (\ref{m3}) to get for some $b<\frac{1}{2}$:
\begin{eqnarray*}
& & \n{(P_{Q^l_{\alpha}}J^s v_i)(P_l \Lambda ^{\frac{1}{2}}v_j) (P_{\Delta l} v_k) \prod_{\nu  \neq i,j,k} P_l v_{\nu}}{\XX{0}{-b}} \\
& \le & c 2^{ls_0} \n{(P_{Q^l_{\alpha}}J^s v_i)(P_l \Lambda ^{\frac{1}{2}}v_j) (P_{\Delta l} v_k) \prod_{\nu  \neq i,j,k} P_l v_{\nu}}{L^{p'_0}_{xt}} \\
& \le & c 2^{ls_0} \n{P_{Q^l_{\alpha}}J^s v_i}{L_t^{r_0}(L_x^{q_0})} \n{P_l \Lambda ^{\frac{1}{2}}v_j}{L_t^{2}(L_x^{q})} \n{P_{\Delta l} v_k}{L_t^{r_1}(L_x^{q_1})} \prod_{\nu  \neq i,j,k} \n{P_l v_{\nu}}{L^p_{xt}}\\
& \le & c 2^{l(s_0 + \sigma_1 + \epsilon -s)} \n{P_{Q^l_{\alpha}}J^s v_i}{\XXm{0}{b}}\n{v_j}{\XXm{s}{\frac{1}{2}}} \prod_{k \neq i,j} \n{v_k}{\xm}\,\,.
\end{eqnarray*}
Squaring the last and summing up over $\alpha$ we arrive at the square of (\ref{m14}).

4. Conclusion: Since $s_0 + s_1-s<0$ as well as $s_0 + \sigma_1 + \epsilon -s<0$ we can now insert (\ref{m13}) and (\ref{m14}) into (\ref{m12}) and finish the proof by summing up over $k$ and $l$.
$\hfill \Box$

\begin{lemma}\label{l33} Let $m,n \in \N$ with $m \ge 2$, $m+n \ge 4$ and $s> \frac{n}{2} -\frac{1}{m-1}$. For $1 \le i,j \le m$ and $v_i \in \XXm{s}{\frac{1}{2}}$ define $f_i(\xi, \tau)= <\xi>^s <\tau - |\xi|^2>^{\frac{1}{2}}\F v_i (\xi, \tau)$ and
\[G_{0j}(\xi,\tau)=<\tau + |\xi|^2>^{-\frac{1}{2}} \int d \nu <\xi_j>^s \chi_A \prod_{i=1}^m<\tau_i - |\xi_i|^2>^{-\frac{1}{2}}<\xi_i>^{-s}f_i(\xi_i,\tau_i)\,\,,\]
where in $A$ the inequality $<\tau + |\xi|^2> \ge \max_{i=1}^m <\tau_i - |\xi_i|^2>$ holds. Then there exists a $b<\frac{1}{2}$ for which the estimate
\[\n{G_{0j}}{L^2_{\xi}(L^1_{\tau})} \le c \prod_{i=1}^m \n{v_i}{\xm}\]
is valid.
\end{lemma}

Proof: We choose $\epsilon \in (0,s-\frac{n}{2} +\frac{1}{m-1})$ with $\epsilon \le \frac{1}{m-1}$ and define $\delta=\frac{m-1}{2m}\epsilon$. Observe that, because of
\[\sum_{i=1}^m <\xi_i>^2 \le <\tau + |\xi|^2> + \sum_{i=1}^m <\tau_i - |\xi_i|^2>\]
in the region $A$ the inequality
\[<\tau + |\xi|^2> \ge c \prod_{i=1}^m <\tau_i - |\xi_i|^2>^{2 \delta}\prod_{i=1 \atop i\neq j}^m <\xi_i>^{\frac{2}{m-1}-2\epsilon}\]
holds. From this we obtain
\[G_{0j}(\xi,\tau) \le c \int d \nu \prod_{i=1 \atop i\neq j}^m <\xi_i>^{-s-\frac{1}{m-1}+\epsilon}\prod_{i=1}^m <\tau_i - |\xi_i|^2>^{-\frac{1}{2} - \delta}f_i(\xi_i,\tau_i)\,\,.\]
In order to estimate $\n{G_{0j}}{L^2_{\xi}(L^1_{\tau})}$ by duality let $f_0 \in L^2_{\xi}$ with $f_0 \ge 0$ and $\n{f_0}{L^2_{\xi}}=1$. By Fubini and Cauchy-Schwarz we get:
\begin{eqnarray*}
& & \int \mu (d \xi ) d \tau d \nu f_0(\xi) G_{0j}(\xi,\tau) \\
& \le & c\int \mu (d \xi ) d \tau d \nu f_0(\xi) \prod_{i=1}^m <\tau_i - |\xi_i|^2>^{-\frac{1}{2} - \delta}f_i(\xi_i,\tau_i) \prod_{i=1 \atop i\neq j}^m <\xi_i>^{-s-\frac{1}{m-1}+\epsilon} \\
& = & \!\!\!\!c\!\! \int\!\! \mu (d \xi_1 .. d \xi_m)d \tau_1 .. d \tau_m f_0(\sum_{i=1}^m\xi_i)\!\!\prod_{i=1}^m\!\! <\!\!\tau_i - |\xi_i|^2\!\!>^{-\frac{1}{2} - \delta}\!\!f_i(\xi_i,\tau_i)\!\! \prod_{i=1 \atop i\neq j}^m \!\!<\!\!\xi_i\!\!>^{-s-\frac{1}{m-1}+\epsilon} \\
& \le & \!\!\!\!c\!\! \int\!\! \mu (d \xi_1 .. d \xi_m)f_0(\sum_{i=1}^m\xi_i)\!\!\prod_{i=1 \atop i\neq j}^m \!\!<\!\!\xi_i\!\!>^{-s-\frac{1}{m-1}+\epsilon} \!\!\prod_{i=1}^m (\int \!\!d \tau_i f_i(\xi_i,\tau_i)^2\!\! <\!\!\tau_i - |\xi_i|^2\!\!>^{-\delta})^{\frac{1}{2}} \\
& \le & c \prod_{i=1 \atop i\neq j}^m (\int \mu (d\xi_i)<\xi_i>^{-2s-\frac{2}{m-1}+2\epsilon})^{\frac{1}{2}}\prod_{i=1}^m \n{f_i <\tau - |\xi|^2>^{-\frac{\delta}{2}}}{L^2_{\xi \tau}} \\
& \le & c \prod_{i=1}^m \n{f_i <\tau - |\xi|^2>^{-\frac{\delta}{2}}}{L^2_{\xi \tau}} = c \prod_{i=1}^m \n{v_i}{\XXm{s}{\frac{1- \delta}{2}}}\,\,.
\end{eqnarray*}
From this the statement of the lemma follows for $b=\frac{1- \delta}{2}$.
$\hfill \Box$

\vspace{0,5cm}

Proof of Theorem \ref{t2}: 1. Setting $v_i = \overline{u_i}$ the claimed estimates read
\begin{eqnarray}
\n{ \prod_{i=1}^m v_i}{\XX{s+1}{-\frac{1}{2}}} \leq c T^{\theta} \prod_{i=1}^m \n{v_i}{\XXm{s}{\frac{1}{2}}} \,\,, \\
\n{ \prod_{i=1}^m v_i}{\yy{s+1}} \leq c T^{\theta} \prod_{i=1}^m \n{v_i}{\XXm{s}{\frac{1}{2}}}\,\,.
\end{eqnarray}
To prove these, we shall assume $s \in (\frac{n}{2}-\frac{1}{m-1},\frac{n}{2})$ as well as $s \le \frac{1}{3}$ for $n=1$ and $m \in \{3,4\}$. Now for $f_i(\xi,\tau)=<\tau - |\xi|^2>^{\frac{1}{2}}<\xi>^s \F v_i (\xi,\tau)$  we have by Lemma \ref{l11}, that the left hand side of (19) is equal to
\begin{eqnarray*}
\n{<\tau + |\xi|^2>^{-\frac{1}{2}} <\xi>^{s+1}\int d \nu  \prod_{i=1}^m <\tau_i - |\xi_i|^2>^{-\frac{1}{2}}<\xi_i>^{-s}f_i(\xi_i,\tau_i)}{L^2_{\xi\tau}}\\
\le c \sum_{i=0}^m\n{F_i}{L^2_{\xi\tau}} \,\,,\hspace{5cm}
\end{eqnarray*}
where
\[F_0(\xi,\tau)=<\xi>^{s}\int d \nu  \prod_{i=1}^m<\tau_i - |\xi_i|^2>^{-\frac{1}{2}}<\xi_i>^{-s}f_i(\xi_i,\tau_i)\]
and for $1 \le i \le m$
\[F_i(\xi,\tau)=<\!\!\tau + |\xi|^2\!\!>^{-\frac{1}{2}}<\!\!\xi\!\!>^{s}\!\!\int \!\! d \nu  <\!\!\tau_i - |\xi_i|^2\!\!>^{\frac{1}{2}}\prod_{k=1}^m<\!\!\tau_k - |\xi_k|^2\!\!>^{-\frac{1}{2}}<\!\!\xi_k\!\!>^{-s}f_k(\xi_k,\tau_k)\,\,.\]
Here we have used the inequality
\[<\xi>^{2} \le <\tau + |\xi|^2> + \sum_{i=1}^m <\tau_i - |\xi_i|^2>\,\,.\]
Now by $<\xi> \le \sum_{j=1}^m <\xi_j>$ it follows, that
\[F_0(\xi,\tau) \le \sum_{j=1}^m F_{0j}(\xi,\tau), \hspace{2cm}F_i(\xi,\tau) \le \sum_{j=1}^m F_{ij}(\xi,\tau)\,\,,\]
where
\[F_{0j}(\xi,\tau) = \int d \nu <\xi_j>^{s} \prod_{i=1}^m<\tau_i - |\xi_i|^2>^{-\frac{1}{2}}<\xi_i>^{-s}f_i(\xi_i,\tau_i)\]
and
\[F_{ij}(\xi,\tau)=<\!\!\tau + |\xi|^2\!\!>^{-\frac{1}{2}}\!\!\int \!\! d \nu  <\!\!\tau_i - |\xi_i|^2\!\!>^{\frac{1}{2}}<\!\!\xi_j\!\!>^{s}\prod_{k=1}^m<\!\!\tau_k - |\xi_k|^2\!\!>^{-\frac{1}{2}}<\!\!\xi_k\!\!>^{-s}f_k(\xi_k,\tau_k)\,\,.\]

2. To derive the estimate (20) we use the inequality
\[<\xi>^{2} \le c(<\tau + |\xi|^2> \chi_A + \sum_{i=1}^m <\tau_i - |\xi_i|^2>)\,\,,\]
where in the region $A$ we have $<\tau + |\xi|^2> \ge \max_{i=1}^m <\tau_i - |\xi_i|^2>$ (cf. Lemma \ref{l33}). Now again by Lemma \ref{l11} we see that the left hand side of (20) is equal to
\begin{eqnarray*}
\n{<\tau + |\xi|^2>^{-1} <\xi>^{s+1}\int d \nu  \prod_{i=1}^m <\tau_i - |\xi_i|^2>^{-\frac{1}{2}}<\xi_i>^{-s}f_i(\xi_i,\tau_i)}{L^2_{\xi}(L^1_{\tau})}\\
\le c \sum_{i=0}^m\n{G_i}{L^2_{\xi}(L^1_{\tau})}\,\,, \hspace{5cm}
\end{eqnarray*}
where now
\begin{eqnarray*}
G_0(\xi,\tau)&=&<\!\!\tau + |\xi|^2\!\!>^{-\frac{1}{2}}<\xi>^{s}\int d \nu  \chi_A \prod_{i=1}^m<\tau_i - |\xi_i|^2>^{-\frac{1}{2}}<\xi_i>^{-s}f_i(\xi_i,\tau_i)\\
& \le & \sum_{j=1}^m G_{0j}(\xi,\tau)
\end{eqnarray*}
with $G_{0j}$ precisely as in Lemma \ref{l33}, and for $1 \le i \le m$
\[G_i(\xi,\tau)=<\!\!\tau + |\xi|^2\!\!>^{-1}<\!\!\xi\!\!>^{s}\!\!\int \!\! d \nu  <\!\!\tau_i - |\xi_i|^2\!\!>^{\frac{1}{2}}\prod_{k=1}^m<\!\!\tau_k - |\xi_k|^2\!\!>^{-\frac{1}{2}}<\!\!\xi_k\!\!>^{-s}f_k(\xi_k,\tau_k)\,\,.\]
Using Cauchy-Schwarz' inequality the estimation of $G_i$, $1 \le i \le m$, can easily be reduced to the estimation of $F_i$, in fact for any $\epsilon > 0$ we have:
\[\n{G_i}{L^2_{\xi}(L^1_{\tau})} \le c_{\epsilon}\n{<\tau + |\xi|^2>^{\epsilon}F_i}{L^2_{\xi\tau}}\le \sum_{j=1}^m c_{\epsilon}\n{<\tau + |\xi|^2>^{\epsilon}F_{ij}}{L^2_{\xi\tau}}\] 

3. Using Lemma \ref{l11} from the introduction and Lemma \ref{l31} we have for \\ $1 \le j \le m$:
\[\n{F_{0j}}{L^2_{\xi\tau}}=c\n{(J^s v_j)\prod^m_{i=1, i \neq j} v_i}{L^2_{xt}} \leq c \prod_{i=1}^m \n{v_i}{\xm}\]
for some $b<\frac{1}{2}$. Now we use Lemma \ref{l12} to conclude that
\[\n{F_{0j}}{L^2_{\xi\tau}}\leq c T^{\theta}\prod_{i=1}^m \n{v_i}{\XXm{s}{\frac{1}{2}}}\]
for some $\theta > 0$. Similarly, but using Corollary \ref{k31} (resp. Lemma \ref{l32}) instead of Lemma \ref{l31}, we get the same upper bound for $\n{<\tau + |\xi|^2>^{\epsilon}F_{ij}}{L^2_{\xi\tau}}$, provided $\epsilon$ is sufficiently small, for $1 \le i = j \le m$ (resp. $1 \le i \neq j \le m$). Now the estimate (19) is proved. For the proof of (20), it remains to show that $\n{G_{0j}}{L^2_{\xi}(L^1_{\tau})}$, $1 \le j \le m$, is bounded by the same quantity. But this follows by Lemma \ref{l33} and Lemma \ref{l12}.
$\hfill \Box$


\begin{thebibliography}{99}

\bibitem[B93]{B93}Bourgain, J.: Fourier transform restriction phenomena for certain lattice subsets and applications to nonlinear evolution equations, GAFA 3 (1993), 107 - 156 and 209 - 262
\bibitem[G96]{G96}Ginibre, J.: Le probl\`eme de Cauchy pour des EDP semi-lin\'eaires p\'eriodiques en variables d'espace (d'apr\`es Bourgain), Ast\'erisque 237 (1996), 163 - 187
\bibitem[GTV97]{GTV97}Ginibre, J., Tsutsumi, Y., Velo, G.: On the Cauchy-Problem for the Zakharov-System, J. of Functional Analysis 151 (1997), 384 - 436
\bibitem[KPV96]{KPV96}Kenig, C. E., Ponce, G., Vega, L.: Quadratic forms for the 1 - D semilinear Schr\"odinger equation, Transactions of the AMS 348 (1996), 3323 - 3353
\bibitem[T99]{T}Takaoka, H.: Well-posedness for the one dimensional nonlinear Schr\"odinger equation with the derivative nonlinearity, Advances in Differential Equations 4 (1999), 561 - 580
\bibitem[HW]{HW}Hardy, G. H.; Wright, E. M.: Einf\"uhrung in die Zahlentheorie; M\"unchen, 1958
\bibitem[P]{P}Petersson, H.: Modulfunktionen und quadratische Formen; Berlin - Heidelberg - New York, 1982
\end{thebibliography}
\end{document}